\documentclass{amsart}

\usepackage{amssymb, amscd, amsmath, amsthm, epsf, epsfig, latexsym}

\def\zz{{\bf Z}}
\def\qq{{\bf Q}}

\def\rr{{\bf R}}
\def\nsmall{ \mbox{\tiny{N} }}

\newtheorem{theorem}{Theorem}[section]
\newtheorem{lemma}[theorem]{Lemma}

\newtheorem{conjecture}[theorem]{Conjecture}

\theoremstyle{definition}

\newtheorem{example}[theorem]{Example}

\theoremstyle{remark}

\def\dfc{\mbox{\rm def}}
\begin{document}

\title{Four-Manifolds of Large Negative Deficiency}

\author{Charles Livingston}

\address{Department of Mathematics, Indiana University, Bloomington, IN
47405}

\email{livingst@indiana.edu}

\keywords{Homology sphere, deficiency,  contractible manifold}

\subjclass{}

\date{\today}

 \begin{abstract}For every   $N > 0$ there exists a group of deficiency
less than
${-}N$ that arises as the fundamental group of a smooth homology
4--sphere and also as the fundamental group of the complement of a
compact contractible submanifold of the 4--sphere.  A group is the
fundamental group of the complement of a contractible submanifold of the
$n$--sphere,
$n > 4$, if and only if it is the fundamental group of a homology
$n$--sphere. There exist fundamental groups of homology $n$--spheres,
$n > 4$, that cannot arise as the fundamental group of the complement of
a contractible submanifold of the 4--sphere.
 \end{abstract}

\maketitle

\vskip.2in 
\section{Introduction}

The deficiency of a finite presentation of a group $G$   is $g - r$,
where $g$ is the number of generators in the presentation and $r$ is
the number of relations.    The deficiency of $G$, $\mbox{def}(G)$, is
 the maximum value of this difference, taken over all
presentations.     The deficiency of a manifold, $\dfc(X)$, is defined
to be $\dfc(\pi_1(X))$.  Our goal is the proof  
 of the following two theorems.

\begin{theorem}\label{hs} For all $N > 0$ there exists a smooth
4--dimensional  homology sphere
$X_{\nsmall}^4$ such that $\dfc(X_{\nsmall}^4) < {{-}N}$.
\end{theorem}

\begin{theorem}\label{cont} There exists a smooth, compact,
contractible 4--manifold
$Z$ and for all $N > 0$   an embedding $f_{\nsmall} : Z \to S^4$ such
that
$\dfc (S^4 - f_{\nsmall} (Z)) < {{-}N}$. Furthermore, $Z  \times I \cong
B^5$ so, in particular, $Z$ embeds in $S^4$ with contractible
complement.
\end{theorem}

Kervaire~\cite{ke} proved that a group $G$ is the fundamental group  
of a homology sphere  of dimension greater than four if and only if $G$
is finitely presented and
$H_1(G) = 0 = H_2(G)$.   Hausmann and Weinberger~\cite{hw} demonstrated
the existence of groups satisfying these conditions that cannot occur
as the  fundamental group of a homology 4--sphere. These groups have
arbitrarily large negative deficiency  and provided the first examples
of the analog of Theorem~\ref{hs} in high dimensions.  Kervaire also 
proved that all perfect groups with deficiency 0 occur as the
fundamental groups of homology 4--spheres and Hillman~\cite{h,h2} recently
constructed examples of homology 4--spheres of deficiency ${-}1$.

Lickorish~\cite{lk} proved that any finitely presented perfect group
$G$ with
$\dfc(G) = 0$ occurs as the fundamental group of the complement of a
contractible   submanifold of $S^4$.  In \cite{lv}, Lickorish's
construction was modified to build a contractible 4--manifold with an
infinite number of embeddings into $S^4$, distinguished by the
fundamental groups of the complement. The proof of Theorem~\ref{cont}
generalizes the construction of \cite{lv}, giving the first complements
of contractible 4--manifolds with negative deficiency.

Theorem~\ref{hs} occurs naturally in the context of Kervaire's theorem. 
Theorem~\ref{cont} can similarly be placed in the context of higher
dimensional topology with the following theorem. 
 \begin{theorem}\label{hdim} For any $n  \ge 5$ and    group
$G$, there exists a smooth, compact, contractible,  $n$--dimensional
submanifold  $M_G \subset S^n$ with 
$\pi_1(S^n - M_G ) = G$ if and only if $G$ is finitely presented,
perfect, and
$H_2(G) = 0$.  In this case  $M_G \times I
\cong B^{n+1}$, so, in particular, there exists an embedding 
$\phi: M_G \to S^n$  such that $S^n - \phi(M_G) $ is contractible.
\end{theorem}

An extension of Hausmann and Weinberger's work yields the following
theorem.
\begin{theorem}\label{theoremhw} There exist  groups $G$ that occur as
the fundamental groups of the complements of compact contractible
submanifolds in $S^n$  for 
$n \ge 5$ but cannot occur as the fundamental groups of the complements
of compact contractible manifolds in $S^4$.
\end{theorem}

\noindent{\bf Notation and Conventions} We work  in the smooth
category.  Homology and cohomology are taken with integer  coefficients
unless otherwise noted. The ring $ \zz[t, t^{{-}1}]$ is denoted
$\Lambda$.

\section{Deficiency}\label{cover}

Levine  constructed   knots in $S^4$ with complements of arbitrarily
large negative deficiency
\cite{le}. His proof was  based on the example of the 2--twist spin of
the trefoil knot, for which he showed the complementary group has
deficiency
${-}1$.  The construction of 4--manifolds used here begins with a manifold
$W^4$   which is then modified by replacing the tubular neighborhood of
a collection of embedded circles with complements of the 2--twist spin
of the trefoil.  

Levine's analysis of the deficiency of the groups    depended on
computations of the homology of the infinite cyclic cover of the
corresponding spaces.  This approach is not immediately available
here,  since the relevant spaces have trivial first homology and hence
no infinite  cyclic cover.  The  way around this difficulty is to work
with manifolds that have finite covers which themselves do have
infinite cyclic covers. 

This section contains the necessary covering space theory, its
application to the study of deficiency, and a review of Levine's
results concerning the deficiency of
$\Lambda$--modules.  In the next section the applicable covering spaces
of the 2--twist spin of the trefoil are studied.  In
Section~\ref{construct} the general construction is described and the
results of sections~\ref{cover} and \ref{tref} are applied to prove
theorems~\ref{hs} and
\ref{cont}.

\subsection{Deficiency of Finite Index Subgroups}
 
 \begin{theorem}\label{finitethm} If $H$ is a subgroup of index $n$ in
a finitely presented group $G$, then  $\dfc(H) \ge n \dfc(G)  - n +
1$. 
\end{theorem}
\begin{proof} A purely algebraic proof of this is based on the
Reidemeister--Schreier rewriting process.  Details are contained in 
\cite{kms}, or see \cite{st} for a covering space description. 
Briefly, $G$ is the fundamental group of a 2--complex with one vertex,
$g$ 1--cells and $r$ 2--cells.  Hence, the $n$--fold cover with group
$H$ has a $1$--skeleton built from $n$ 0--cells and $ng$ 1--cells, so
is homotopy equivalent to a 1--complex with one vertex and $ng-n+1$
1--cells.  This has free fundamental group with $ng - n+1$ generators. 
There are $nr$ 2--cells, so $H$ has a presentation with $ng - n+1$
generators and $nr$ relations, and hence deficiency $n(g-r) - n +1$.
The  result follows.
\end{proof}

\subsection{Deficiency and Infinite Cyclic Covers}

Suppose that there is a surjective  homomorphism $\phi$ from a finitely
presented group
$G$ onto
$\zz$  with kernel $H$.  There is an action of $\zz$ on the
abelianization of $H$,
$H_1(H)$, making $H_1(H)$ into a $\Lambda$--module.  In terms of
covering space, if
$K_G$ is a space with fundamental group $G$, then $H_1(H)$ is the first
homology of the (connected) infinite cyclic cover of $K_G$
corresponding to $H$ and the
$\Lambda$--action is induced by the group of deck transformations.

\vskip.1in
\noindent{\bf Deficiency of $\Lambda$--modules.} For an arbitrary
finitely presented
$\Lambda$--module $M$ the deficiency is defined as for groups.  The
deficiency of a presentation is the difference of the number of
generators and the number of relations.  The deficiency of the module
is the maximum value of this difference, taken over all presentations. 
It is bounded above by the rank of $M$; that is, by the dimension of the
$\qq(t)$--vector space $M \otimes \qq(t)$.
\vskip.1in

 As just described, if $\phi: G \to \zz$ is a surjective homomorphism,
the abelianization of the kernel, $H$, is a $\Lambda$--module,
$H_1(H)$.  In this situation there is  the following bound.
 
\begin{theorem}\label{infthm}    $\dfc( H_1(H) )
\ge
\dfc(G) -1$.
\end{theorem}
\begin{proof} If $G$ has a presentation with $g$ generators and $r$
relations, then via the Reidemeister-Schreier rewriting process, or the
corresponding geometric construction, 
$H_1(H)$ has a presentation with $g-1$   generators and
$r$ relations as a $\Lambda$--module.    The deficiency of this
presentation is
$g-1  -r =  (g-r )- 1$.  The result follows. 
\end{proof}

\vskip.1in
\subsection{Iterated Covers and Deficiency}

\begin{theorem} \label{bigthm} Suppose   $W$ is a connected manifold, 
$\tilde{W}$ is  a connected
$n$--fold cover of $W$, and  $\phi: \pi_1(\tilde{W}) \to \zz$ is a
surjection.  If
$\tilde{W}_\infty$ is the associated infinite cyclic cover of
$\tilde{W}$,   then $$
\dfc(H_1(\tilde{W}_\infty) )\ge  n(\dfc(\pi_1(W)) - n.   
$$
\end{theorem}

\begin{proof} The result is an immediate corollary of
theorems~\ref{finitethm}~and~\ref{infthm}.
\end{proof}

\subsection{The Deficiency of $\Lambda$--Modules}

For a finitely presented $\Lambda$--module  the following inequality
was proved by Levine.  The proof is summarized below, with relevant
definitions embedded in the  summary.  Details can be found in
\cite{le}.

\begin{theorem}\label{le}If  $M$ a finitely presented $\Lambda$--module
of rank $r_0$ and deficiency  $d$, then  $\mbox{\rm
Ext}_\Lambda^2(M,\Lambda)$ can be generated, as a
$\Lambda$--module, by $r_0 - d$ elements.

\end{theorem}
\begin{proof} A finite presentation of $M$ with $g$ generators and $r$
relations ($g - r = d$)   yields an  exact sequence
$$0
\to
 {F_2}
\to  F_1
\xrightarrow{\phi_1} F_0 \to M \to 0,$$ with $F_0$ free of rank $g$ and
$F_1$ free of rank
$r$. By definition, $F_2$ is the kernel of $\phi_1$  and homological
properties of
$\Lambda$  imply that
$F_2$ is free. Tensoring with
$\qq(t)$ yields an exact sequence of $\qq(t)$ vector spaces with
dimensions, by definition, the ranks of each of the modules.  Since the
alternating sum of these dimensions is 0,   rank($F_2$) = $r - g + r_0
= r_0 - d$.

By definition, $\mbox{Ext}_\Lambda^2(M,\Lambda)$ is the cokernel of the
induced map
$\mbox{Hom}_\Lambda(F_1, \Lambda) \to  \mbox{Hom}_\Lambda(F_2,
\Lambda)$, and hence is
  a quotient of a free module of rank $r_0-d$.  It follows that it  can
be generated by
$r_0 - d$ elements, as desired.
\end{proof}

\begin{example} The $\Lambda$ modules $\Lambda$, $\Lambda/\left< 3
\right> $ and
$\Lambda/\left<3  , t + 1\right>$ have deficiencies $1$, $0$, and ${-}1$,
respectively. 
\end{example}

The obvious presentations have these deficiencies, so it remains to see
that these represent the maximum deficiencies.  For $\Lambda$ itself
this is immediately seen by moving to the vector space setting via
tensoring with $\qq(t)$.  For $\Lambda/\left< 3 \right> $ the same
argument applies, since  $\Lambda/\left< 3 \right> \otimes
\qq(t) = 0$.

For $\Lambda/\left<3\ ,\ t + 1\right>$, notice that its rank, $r_0$, is
0.  Hence, by applying Theorem~\ref{le} the question is reduced to
determining whether
$\mbox{Ext}_\Lambda^2(\Lambda/\left<3  ,  t + 1\right>,\Lambda)$ can be
generated by fewer than 1 element; that is, is
$\mbox{Ext}_\Lambda^2(\Lambda/\left<3  ,  t + 1\right>,\Lambda)$
trivial?  The free resolution of
$ \Lambda/\left<3  ,  t + 1\right>$ is given by 
$$0
\to
\Lambda
\xrightarrow{\phi_2} \Lambda^2
\xrightarrow{\phi_1} \Lambda \to  \Lambda/\left<3  ,  t + 1\right> \to
0,$$ where
$\phi_1( (1,0)) = 3$, $\phi_1( (0,1)) = t+1$, and $\phi_2( 1) = (t+1,
{-}3)$.  Applying Hom  gives that $\mbox{Ext}_\Lambda^2(\Lambda/\left<3 
,  t + 1\right>,\Lambda)$ is the cokernel of the map $$\Lambda^2
\xrightarrow{\phi_2^*} \Lambda,$$ where $\phi_2^*((1,0) = t+1$ and
$\phi_2^*((0,1) = {-}3$.  Hence, $\mbox{Ext}_\Lambda^2(\Lambda/\left<3 
,  t + 1\right>,\Lambda) \cong 
\Lambda/\left< t+1   , {-}3\right>$. As desired, this is nontrivial---as
an abelian group it is isomorphic to $\zz_3$.

\section{Building blocks for high deficiency examples}\label{tref}

Let $K$ denote the 2--twist spin of the trefoil knot in $S^4$.  Let
$E$ denote the complement of an open tubular neighborhood of $K$.
This space is a fiber bundle over
$S^1$ with fiber a punctured lens space, $L(3,1)_0$. The monodromy is
an involution which acts nontrivially on $\zz_3 = H_1(L(3,1)_0)$ and the
computation of the fundamental group follows readily:

\begin{theorem} $\pi_1(E) = \langle t, x\ |\ x^3 = 1, txt^{{-}1} =
x^{{-}1}
\rangle$.
\end{theorem}

There is a homomorphism of $\pi_1(E)$ to $\zz$ sending $t$ to $d$ and
$x$ to 0.  This induces an infinite cyclic cover of $E$,   denoted
$\tilde{E}_{d}$.  The usual example is the case of $d = 1$, for which
$\tilde{E}_{d} \cong  L(3,1)_0 \times \rr$.  If
$d = 0$ the cover consists of an infinite family of copies $E$.  If
$d > 1$ the cover is disconnected, consisting of $d$ copies of
$\tilde{E}_{1}$.

Computing the first homology of this cover as a $\Lambda$--module
 follows from standard methods in the case of
$d = 1$ and is a trivial calculation in the case of $d=0$.  For $d >1$
the deck transformation $T$ cyclically permutes the $d$ copies of
$\tilde{E}_{1}$, with $T^d$ restricting to give the generating deck
transformation on each of the copies of
$\tilde{E}_{1}$.  The next theorem is a consequence of these
observations.

 \begin{theorem} For $d > 0$, $H_1(\tilde{E}_{d}) \cong \Lambda /
\langle 3 , t^d +1 \rangle$.  For $d = 0$, 
 $H_1(\tilde{E}_{d})
\cong \Lambda $.   
\end{theorem}

In order to apply a Mayer-Vietoris argument it is also necessary to
understand the homology of the boundary of covers of $E$.  However,
the boundary of $E$ is $S^1
\times S^2$, and the $d > 0$ infinite cyclic covers of the boundary
have trivial first homology.  In general:

\begin{theorem} For $d >0$, $H_1(\partial \tilde{E}_{d} ) \cong 0$. 
For
$d = 0$, $H_1(\partial \tilde{E}_{d} ) \cong \Lambda$.

\end{theorem}

Let $W$ be  a compact 4--manifold and $\phi: H_1(W) \to \zz$ be a
surjection inducing an infinite cyclic cover, $\tilde{W}$.  Suppose that
$\alpha$ is an oriented simple closed curve in $W$ and $\phi([\alpha]) =
d$.   Let 
$W^*$ be the manifold constructed by removing a tubular neighborhood of
$\alpha$ and replacing it with
$E$ via any homeomorphism of the boundaries.  Removing
$\alpha$ does not change the first homology of $W$ or of $\tilde{W}$.  A
Mayer-Vietoris argument gives the homology of the infinite cyclic cover
$\tilde{W}^*$.  

\begin{theorem} If  $d >0$, $H_1(\tilde{W}^*) \cong  H_1(\tilde{W} )
\oplus\Lambda/\langle 3, t^d +1 \rangle$.  If $d = 0$, 
$H_1(\tilde{W}^*) \cong  H_1(\tilde{W} )$.   
\end{theorem}


\section{Building the Examples}\label{construct}

The starting point for our construction is a perfect group with
specified properties.  The simplest example is now described.  Let
$G$  
 be the free product $D * D$ where $D$ is the binary icosahedral group,
the fundamental group of  the Poincar\'e dodecahedral homology
3--sphere. The group $D$ has    order 120, is perfect, and has
deficiency 0. (One presentation is $\langle x,y\ |\ 
 x^2 = y^3 = (xy)^5\rangle$.)  The homomorphism of
$D * D$ to
$D$ given by isomorphisms on each factor has an index 120 kernel
isomorphic to the free product of 119 copies of
$\zz$.  (To see this, let $M$ be the dodecahedral space with the
interior of a 3--ball removed and consider the restriction to $\partial(M
\times I)$ of the (120--fold) universal cover of $M \times I$. It is
clear that this cover corresponds to the kernel of the map on fundamental
groups induced by the inclusion of $\partial(M
\times I)$ into $M \times I$.  Since the universal cover of $M$ is a 120
times punctured 3--sphere, the boundary of the cover of $M \times I$
consists of two copies of $S^3$ each with 120 balls removed and with
boundaries identified.  The resulting space is $\#_{119} S^1 \times S^2$.)

For the proof of Theorem~\ref{hs}, let $W$ be a homology 4--sphere with
fundamental group $G$.  	It exists by \cite{ke}. For the proof of
Theorem~\ref{cont}, let 
$W$ be the complement of a contractible manifold $Z$ embedded in
$S^4$ with $\pi_1(W) = G$.
  By~\cite{lk} such a $Z$ exists, with the additional property that $Z
\times I = B^5$.

 Let $\tilde{W}$ be the associated  120--fold cover of $W$ and let
$\phi:\pi_1(\tilde{W}) \to \zz$ be an arbitrarily chosen surjective
homomorphism.  Let $\tilde{W}_\infty$ be the corresponding infinite
cyclic cover.  

Select an oriented simple closed curve $x$ in $\tilde{W}$ such that
$\phi([x]) = 1$.  It can be assumed via general position that the
projection of $x$,
$\overline{x}$, in $W$ is an embedding.  

Replace each tubular neighborhood of   $k$ parallel copies of
$\overline{x}$ with $E$, the complement of the 2--twist spin of the
trefoil.  Call the resulting manifold
$Y_k$.  In the case that $W$ is a homology 4--sphere, $Y_k$ will also
be a homology 4--sphere. In the case that $W$ is the complement of a
contractible manifold $Z$ in $S^4$, since the
operation of replacing the neighborhood of a circle in
$S^4$ with a knot complement gives $S^4$ back again if the correct
glueing homeomorphism is used, it follows  that
$Y_k$ embeds in
$S^4$ with complement diffeomorphic to $Z$.

Let $\tilde{Y}_{k}$ denote the 120--fold cover of $Y_k$ and let
$\tilde{Y}_{k,\infty} $ denote the infinite cyclic cover of
$\tilde{Y}_{k}$.

 Suppose that
$\pi_1(Y_k)$ has deficiency
$d$.  Then according Theorem~\ref{bigthm} the deficiency of
$H_1(\tilde{Y}_{k,\infty})$ satisfies
$$\mbox{def}(H_1(\tilde{Y}_{k,\infty})) \ge 120(d-1).$$

The space $\tilde{Y_k}$ is built from $\tilde{W}$ by removing parallel
copies of neighborhoods of $x$ and replacing them with copies of $E$.  (On each of these
lifts $\phi$ takes value 1.)  In addition, the other lifts of
$\overline{x}$ are removed and replaced by $E$ also.   (Notice that since
$\tilde{W}$ is a regular cover, the projection is a homeomorphism when
restricted to each lift.)  On each of these lifts the representation
$\phi$ takes various values.

From this we see that as a $\Lambda$--module there is the following
decomposition:
$$H_1(\tilde{Y}_{k,\infty}) =  H_1(\tilde{W}_\infty) \oplus
\left( \frac{\Lambda}{\langle 3, t + 1 \rangle} \right) ^k \oplus \left(
   \oplus_i   
\frac{\Lambda}{\langle 3, t^{d_i} +1  \rangle} \right) ^k
\oplus \left( \oplus_j  
\frac{\Lambda}{\langle 3 \rangle} \right) ^k.$$ (The range of the
indices $i$ and $j$ have not been computed explicitly  and may be
empty.  They will soon drop out of the calculations.)

Suppose now that $H_1(\tilde{Y}_{k,\infty})$ has deficiency $D$ and that
$H_1(\tilde{W}_\infty)$ is generated by $N$ elements.  Then, using the
deficiency
$D$ presentation of $H_1(\tilde{Y}_{k,\infty})$ and adding $N$
relations gives a deficiency $D - N$ presentation of the direct sum 
$$M = \left( \frac{\Lambda}{\langle 3, t + 1 \rangle} \right) ^k \oplus
\left(
\oplus 
\frac{\Lambda}{\langle 3, t^{d_i} +1  \rangle} \right) ^k
\oplus \left( \oplus 
\frac{\Lambda}{\langle 3 \rangle} \right) ^k.$$ Hence, the deficiency
of $M$ satisfies def$(M) \ge D - N$.

The  $\Lambda$--module $M$ is   of rank 0 by definition; it is
$\Lambda$--torsion.  According to Theorem~\ref{le}, in this case 
Ext$_\Lambda^2(M,\Lambda)$ can be generated by ${-}\dfc(M)$ elements. In
particular, its quotient,  $$\mbox{Ext}_\Lambda^2
\left( \left(
\frac{\Lambda}{\langle 3, t + 1 \rangle} \right) ^k  \right) \cong 
\left(
\frac{\Lambda}{\langle 3, t + 1 \rangle} \right) ^k     $$ could be
generated by
${-}\dfc(M)$ elements.    Hence,  ${-}\mbox{def}(M) \ge k$, or
$\mbox{def}(M) \le {-}k$. We now have that $D - N \le {-}k$, or $D \le N
{-}k$.  

Applying Theorem~\ref{bigthm} to this situation now gives that
$\mbox{def}(\pi_1(Y_k))$ satisfies
$$120 \dfc(\pi_1(Y_k)) - 120 \le \dfc(H_1(\tilde{Y}_{k,\infty}) ) \le N
- k .$$ We have
$$  \dfc(\pi_1(Y_k)) \le    
\frac{N + 120 -k}{120} .$$ Since this goes to ${-}\infty$ as $k$ goes to
infinity, the manifolds $Y_k$ have the desired negative deficiencies
and theorems~\ref{hs}~and~\ref{cont} are proved.


\section{Higher Dimensional Complements}

In this section we prove Theorem~\ref{hdim}.
\vskip.1in

\noindent{\bf Theorem~\ref{hdim}.} {\sl 
 For any $n \ge 5$ and   group $G$ there exists a smooth,  compact,
contractible, 
$n$--dimensional submanifold  $M_G \subset S^n$ with  $\pi_1(S^n - M_G
) = G$ if and only if $G$ is finitely presented, perfect, and $H_2(G) =
0$.  In this case  $M_G \times I
\cong B^{n+1}$, and in particular there exists an embedding
$\phi: M_G \to S^n$  such that $S^n - \phi(M_G) $ is contractible.}
\vskip.1in

\begin{proof} The second statement is automatic: if $M_G$ is
contractible then $M_G \times I$ is   a  smooth contractible manifold
with boundary a homotopy sphere.  By the $h$--cobordism theorem the
boundary is diffeomorphic to
$S^{n}$ and contains $M_G$.

One direction of the first statement is immediate.  Given the existence
of such an $M_G$, the complement has finitely generated fundamental
group and is a homology ball.  Hence $G$ is perfect.  Since a
$K(G,1)$ can be built from $S^n - M_G$ by adding cells of dimension 3
and higher,
$H_2(G)$ is a quotient of $H_2(S^n - M_G) = 0 $.

For the proof of the reverse implication, let $N$ be an
$n$--dimensional homology sphere with $\pi_1(N) = G$, the existence of
which is given by \cite{ke}.  Fix a handlebody structure on $N$ and let
$N_k$ be the union of all handles of dimension $k$ and less.

The homology group $H_2(N_2)$ is free, say on $a$ generators, and hence
there are $b$ 3--handles, with $b \ge a$.  The corresponding
presentation matrix  of $H_2(N_3)$ is given by a matrix $A$ with $b$
rows and $a$ columns.  Repeated row operations, corresponding to handle
slides on the 3--handles, yields a matrix $A'$  in upper triangular
form; that is, $A'(i,j) = 0$ if $i > j$.  Since $H_2(N_3) = 0$, the
diagonal entries of this matrix must all be $\pm 1$.  

Let $N'_3$ denote the union of $N_2$ with the first $a$ of the  
3--handles in this new handlebody decomposition.  From the
construction, $\pi_1(N'_3) = G$ and
$N'_3$ is a homology ball.  

In the argument that follows it is necessary   that $N'_3$ be stably
parallelizable; in fact, $N'_3$ has trivial tangent bundle, as seen as
follows. Since $N$ is orientable, the tangent bundle of
$N_1$ is trivial.  The obstruction to trivializing the tangent bundle
over the 2--skeleton of
$N$, or over
$N_2$, is the second Stiefel-Whitney class, $w_2 \in H^2(N, \zz_2) =0$.
The obstruction to extending this trivialization over the 3--skeleton
is in $H^3(N,
\pi_2(SO(n)))$.  But
$\pi_2$ of any Lie group is trivial, so this obstruction vanishes
also.

Since $N'_3$ is a stably parallelizable homology ball of dimension 5 or
higher, basic surgery theory \cite{km} implies that  surgery can be
performed on $N'_3$ to yield a
 contractible manifold $M_G$ without changing the boundary.  The
boundary union $N'_3
\cup_\partial M_G$ is a homotopy sphere, since  the map $\pi_1(\partial
N'_3) \to \pi_1(N'_3)$ is surjective.  (To see this surjectivity, note
that  $N'_3$ is built with handles of dimension 3 and less; hence there
is a dual handlebody structure building    
$N'_3$   from $\partial N'_3 \times I$ using handles of dimension $n-3$
and higher. The surjectivity now follows from the assumption   that
$n-3 \ge 2$.) Denote this homotopy sphere, $ N'_3 \cup_\partial M_G $,
by $\Sigma^n$. The $h$--cobordism theorem implies that
$\Sigma
\# {-}\Sigma
\cong S^n$.  

Let $N_G = N'_3 \# {-}\Sigma$.  Then $N_G \cup_\partial M_G \cong S^n$. 
That is, $M_G$ is a contractible manifold embedded in $S^n$ with
complement $N_G$ having fundamental group $G$.  This completes the
proof.

 \end{proof}

\section{Groups that Cannot Occur in Dimension
4}

The goal of this section is to demonstrate the existence of a finitely
presented group $G$ satisfying $H_1(G) = 0 = H_2(G)$ but such that $G$
does not occur as the fundamental group of a compact homology 4--ball. 
The construction is a simple modification of the work of
\cite{hw}, to which the reader is referred for background.  

Fix a prime $p$ and for a group or space $X$ let  $\beta_i^p(X) $
denote the dimension of
$H_i(X,
\zz_p)$.

\begin{lemma}\label{lemmahw} If $G$ is the fundamental group of compact
4--manifold $M$ with nonempty connected boundary and 
$\chi(M) = 1$ (eg. a rational homology 4--ball) and
$H
\subset G$ is  index $k$, then $$2 +
\beta^p_2(H) - 2\beta^p_1(H) \le 2k.$$
\end{lemma}

\begin{proof} Let $\tilde{M}$ be the $k$--fold cover of $M$
corresponding to
$H$.  Then
$$\chi(\tilde{M}) = 1 - \beta^p_1(\tilde{M} )+ \beta^p_2(\tilde{M} )
-\beta^p_3(\tilde{M} ) = k.$$ By duality, $\beta^p_3( \tilde{M} ) =
\beta^p_1( \tilde{M}, \partial
\tilde{M})$.  From  the exact sequence,
$$H_1(  \partial \tilde{M}, \zz_p) \to H_1(  \tilde{M}, \zz_p )
\xrightarrow{\phi} H_1( 
\tilde{M},
\partial\tilde{M} , \zz_p)
\to H_0(\partial \tilde{M}, \zz_p) \to H_0( \tilde{M} , \zz_p),$$ it
follows that  
$$\beta^p_1(\tilde{M}, \partial \tilde{M}) = (\beta^p_0(\partial
\tilde{M}) - 1) +
\mbox{rank}( \mbox{Im} (\phi)).$$ This in turn implies that
$$\beta^p_1(\tilde{M}, \partial \tilde{M})) \le (k - 1) +
\beta^p_1(\tilde{M}).$$ Combining these shows that
$$\beta^p_3(\tilde{M}) \le (k-1) +
\beta^p_1(H).$$  Hence, $$1 - \beta^p_1(H) + \beta^p_2(H) - ((k-1) +
\beta^p_1(H)) \le k.$$ Simplifying gives the desired inequality.

\end{proof}

In \cite{hw} there is constructed, for each positive integer $k$, a
group $G$ satisfying the conditions of Theorem~\ref{hdim} that contains
a subgroup $A$ of finite index, where that index is independent of
$k$.  These groups have the property that $\beta^p_1(A)$ is bounded
above by a linear function of
$k$ and
$\beta^p_2(A)$ is bounded below by a quadratic function of $k$.  
 The inequality of Lemma~\ref{lemmahw} cannot hold for all $k$, and
hence these examples yield the proof of Theorem~\ref{theoremhw}. 

\section{Open problems}

There remains no characterization of the fundamental groups of homology
4--spheres.  In addition, it is unknown whether such a characterization
must depend on the category, smooth or topological.  

Since in higher dimensions the set of fundamental groups of homology
spheres is the same as the set of fundamental groups of complements of
contractible manifolds, one can wonder if the same holds in dimension
4.  We can state this as a conjecture.

\begin{conjecture}A group   is the fundamental group of a homology
4--sphere if and only  if it is the fundamental group of the complement
of a contractible submanifold of the 4--sphere.
\end{conjecture}

Neither direction of the conjecture is known to be true.

\newcommand{\etalchar}[1]{$^{#1}$}

\end{document}